\documentclass[12pt,twoside]{article}
\usepackage{fullpage}
\usepackage[utf8]{inputenc}
\usepackage{macros}
\usepackage{paperstyle}

\title{Riemannian 3-spheres that are hard to sweep out by short curves}
\author{Omar Alshawa and Herng Yi Cheng}

\begin{document}

\maketitle

\begin{abstract}
    We construct a family of Riemannian 3-spheres that cannot be ``swept out'' by short closed curves. More precisely, for each $L > 0$ we construct a Riemannian 3-sphere $M$ with diameter and volume less than 1, so that every 2-parameter family of closed curves in $M$ that satisfies certain topological conditions must contain a curve that is longer than $L$. This obstructs certain min-max approaches to bound the length of the shortest closed geodesic in Riemannian 3-spheres.
    
    We also find obstructions to min-max estimates of the lengths of orthogonal geodesic chords, which are geodesics in a manifold that meet a given submanifold orthogonally at their endpoints. Specifically, for each $L > 0$, we construct Riemannian 3-spheres with diameter and volume less than 1 such that certain orthogonal geodesic chords that arise from min-max methods must have length greater than $L$.
\end{abstract}

\section{Introduction} \label{Introduction}

How long is the shortest closed geodesic in a closed Riemannian manifold $M$ that is not a constant loop? M.~Gromov asked whether its length, denoted by $\scg(M)$, can be bounded in terms of volume $\vol(M)$ of $M$ \cite{Gromov_FillingRiemannianManifolds}. If that is impossible, can $\scg(M)$ be bounded by some function of $\vol(M)$ and the diameter $\diam(M)$ of $M$?

Consider the example where $M$ is a Riemannian 3-sphere. Suppose that $M$ can be ``swept out'' by short curves, in the following sense. The unit 3-sphere $S^3 \subset \R^4$ is partitioned into subspaces $S^1_{st} = S^3 \cap (\R^2 \times \{(s,t)\})$ for $s, t \in \R$ which are circles when $s^2 + t^2 < 1$ and points when $s^2 + t^2 = 1$. Suppose that there exists a continuous map $F : S^3 \to M$ with nonzero degree such that every closed curve $F(S^1_{st})$ is shorter than $\phi(\vol(M), \diam(M))$ for some function $\phi$. Then a standard technique called the \emph{min-max method} would imply that $\scg(M) \leq \phi(\vol(M), \diam(M))$. This method begins from the intuitive picture of $M$ being ``swept out'' by a 2-parameter family of ``short'' closed curves $F(S^1_{st})$. Roughly speaking, a length-shortening process is applied to the entire family of curves $F(S^1_{st})$ continuously, and the property that $\deg F \neq 0$ would force one of the curves to converge to a short but non-constant closed geodesic of length at most $\phi(\vol(M), \diam(M))$.

However, we have ruled out this approach to bounding $\scg(M)$. We accomplished this by constructing Riemannian 3-spheres with small diameter and volume, but which cannot be swept out by short closed curves in the aforementioned manner.

\begin{theorem}[Main result]
    \label{thm:theorem}
    For any $L > 0$, there exists a Riemannian 3-sphere $M = (S^3,g)$ of diameter and volume at most $1$ with the following property: for any continuous map $F : S^3 \to M$ with nonzero degree, one of the closed curves $F(S^1_{st})$ must be longer than $L$.
\end{theorem}

\begin{remark}
    \label{rem:Generalization}
    More general versions of \cref{thm:theorem} can be proven due to the following observation. A crucial ingredient in proving \cref{thm:theorem} is the fact that $\bigcup_s S^1_{st}$ is a 2-sphere for all $0 < t < 1$, because we will apply the Jordan curve theorem to those 2-spheres. In fact, it can be verified that our proof of this theorem will work even when $\{S^1_{st}\}_{s,t}$ is replaced by certain more general foliations $\{Z_{st}\}_{s,t}$ of $S^3$ by 1-cycles $Z_{st}$, as long as for ``almost all'' $t$, $\bigcup_s Z_{st}$ is a 2-sphere.
\end{remark}

\Cref{thm:theorem} should be appraised within the broader context of min-max methods and their applications to the study of the geometry of closed geodesics, minimal submanifolds, and other minimal objects. We will briefly survey these results in this Introduction and suggest potential implications of our theorem. We will also explain how our constructions were motivated by counterexamples to related conjectures.

The ideas behind our proof of \cref{thm:theorem} also helped us to prove a result related to the geometry of \emph{orthogonal geodesic chords}, which may be seen as a ``relative'' analogue of closed geodesics. Given a Riemannian manifold $M$ and a closed submanifold $N \subset M$, an orthogonal geodesic chord is a geodesic in $M$ that starts and ends on $N$, and which is orthogonal to $N$ at its endpoints. What is the length of the shortest orthogonal geodesic chord in this situation? In \cref{sec:Intro_SweepoutsByPaths} we will present another result, \cref{thm:OrthoGeodesicChords}, which obstructs certain min-max approaches to estimating the length of the shortest orthogonal geodesic chord.

When $M$ is a Riemannian manifold that is not simply connected, $\scg(M)$ is at most the \emph{systole} of $M$, which is the length of the shortest non-contractible loop. Gromov proved that when $M$ satisfies a topological condition called \emph{essentialness}, its systole is at most $c(n) \sqrt[n]{\vol(M)}$ where $n = \dim M$ and $c(n)$ is a dimensional constant \cite{Gromov_FillingRiemannianManifolds}. Essential manifolds include all closed surfaces that are not simply connected, as well as all real projective spaces and all tori. A.~Nabutovsky showed that the inequality holds for $c(n) = n$ \cite{Nabutovsky_Systole}.

On the other hand, when $M$ is simply connected, $\scg(M)$ cannot be bounded by a systole. In this sense, it is more difficult to bound $\scg(M)$ when $M$ is a Riemannian sphere. Nevertheless, such bounds are possible for Riemanian 2-spheres $(S^2,g)$. For instance, $\scg(S^2,g) \leq 4\diam(S^2,g)$, as proven by Nabutovsky and R.~Rotman \cite{NabutovskyRotman_ShortestClosedGeodesicS2_Diameter}, and independently by S.~Sabourau \cite{Sabourau_ShortestClosedGeodesicS2_Diameter}. Rotman also proved that $\scg(S^2,g) \leq 4\sqrt2\sqrt{\area(S^2,g)}$ \cite{Rotman_ShortestClosedGeodesicS2_Area}. These results are the sharpest known versions of such bounds for 2-spheres, which were first proven by C.\,B.~Croke and then strengthened by several authors; the early history around these bounds is surveyed in \cite[Section~4]{CrokeKatz_ShortestClosedGeodesic_Survey}.

For manifolds $M$ of dimension at least 3, bounds on $\scg(M)$ have been proven under various conditions on the curvature of $M$ in \cite{BallmanThorbergssonZiller_ShortestClosedGeodesic_Pinched,Tribergs_ShortestClosedGeodesicConvexHypersurface,Croke_ShortestClosedGeodesicConvexHypersurface,NabutovskyRotman_scg_curvature,WuZhifei_ShortestClosedGeodesic4Mfd_Curvature,Rotman_ShortestClosedGeodesic_Curvature,Rademacher_CritVals_Curvature}. On the other hand, no curvature-free bounds for $\scg(M)$ are known for simply-connected manifolds $M$ of dimension at least 3. For such manifolds, the min-max methods sketched before \cref{thm:theorem} are a main avenue for studying $\scg(M)$.

\subsection{Min-max theory for sweepouts by free loops}

The approach to bounding $\scg(M)$ by min-max methods is essentially a quantitative version of the proof that $M$ contains a closed geodesic. The existence of closed geodesics in $M$ was proven by applying a ``geometric calculus of variations'' to the free loop space $\Lambda M$, which is the space of continuous maps $S^1 \to M$ given an appropriate topology \cite{Milnor}. Let $\Lambda^0 M$ denote the space of constant loops on $M$. Since closed geodesics are critical points in $\Lambda M$ with respect to the length functional, we can, in some sense, find those critical points by applying the calculus of variations to continuous maps $f : (D^d, \partial D^d) \to (\Lambda M, \Lambda^0 M)$, where $f$ represents a nontrivial element in $\pi_d(\Lambda M, \Lambda^0 M)$, and $d \geq 1$ is the minimal degree such that $\pi_d(\Lambda M, \Lambda^0 M) \neq 0$. We call such a map $f$ a \emph{sweepout of $M$ by free loops}. The notion of a sweepout allows us to define the min-max value
\begin{equation}
    \label{eq:MinMaxValue_HomotopicalSweepoutFreeLoops}
        \lambda(M) = \inf_{\substack{
            f : (D^d, \partial D^d) \to (\Lambda M, \Lambda^0 M) \\
            0 \neq [f] \in \pi_d(\Lambda M, \Lambda^0 M)
        }} \sup (\length \circ f).
\end{equation}

It can be proven that $\lambda(M)$ is the length of some nonconstant closed geodesic in $M$. Indeed, a similar argument was used by A.\,I.~Fet and L.\,A.~Lyusternik to prove that $M$ must have at least one nonconstant closed geodesic \cite{LyusternikFet_ClosedGeodesicExistence}.\footnote{Surveys of their proof are available in \cite{Bott_MorseTheoryOldNew,Oancea_ClosedGeodesicsExistence}.} Thus $\scg(M) \leq \lambda(M)$, and a natural strategy to bound $\scg(M)$ is to bound $\lambda(M)$ by constructing a single sweepout of $M$ that is ``efficient'' in the sense that every free loop in the sweepout has bounded length.

Such a bound on $\lambda(M)$ may be possible for some combinations of geometric parameters of $M$ but not others. For example, as a consequence of work by Y.~Liokumovich, Nabutovsky, and Rotman, every Riemannian 2-sphere $M = (S^2,g)$ satisfies $\lambda(M) \leq 664\sqrt{\area(S^2,g)} + 2\diam(S^2,g)$ \cite[Theorem~1.3]{liokumovich2015contracting}. Yet $\lambda(S^2,g)$ cannot be bounded in terms of diameter alone or area alone: Sabourau presented a family of Riemannian 2-spheres $(S^2,g)$ for which, after scaling, $\area(S^2,g) \leq 1$ but $\lambda(S^2,g)$ grows without bound \cite[Remark~4.10]{Sabourau_ShortestClosedGeodesicS2_Diameter}. In a similar vein, Liokumovich constructed a family of Riemannian 2-spheres $(S^2,g)$ for which $\diam(S^2,g) \leq 1$ but $\lambda(S^2,g)$ grows without bound \cite{liokumovich2013spheres}. Evidently, the geometry of closed geodesics in $M$ can be probed by either constructing efficient sweepouts, or by constructing ``pathological'' Riemannian manifolds on which any sweepout must be inefficient.

\Cref{thm:theorem} gives a family of Riemannian 3-spheres $M$ for which $\diam(M), \vol(M) \leq 1$ but $\lambda(M)$ grows without bound. Such results proving that $\lambda(M)$ is ``large'' should be understood with the caveat that even if $\lambda(M)$ is large, $M$ may still contain short closed geodesics. This is illustrated by the aforementioned results on $\scg(S^2,g)$ and $\lambda(S^2,g)$. Furthermore, min-max values are in some sense the lengths of closed geodesics of ``positive Morse index,'' so they do not lead to estimates of the lengths of closed geodesics that are local minima of the length functional.

Our result also does not imply that $\scg(M)$ cannot be bounded by a fixed function of $\diam(M)$ and $\vol(M)$. What we can conclude is that if such a bound can be proven using min-max methods, then it is not enough to use sweepouts that represent nontrivial elements of $\pi_2(\Lambda M, \Lambda^0 M)$. More complicated types of sweepouts would be necessary, such as maps $f : (X, X_0) \to (\Lambda M, \Lambda^0 M)$ where $X$ is not a disk, or sweepouts by objects like 1-cycles that generalize free loops.

\subsection{Comparing our constructions to disks whose boundaries are hard to contract}
\label{sec:MotivatingOurConstruction}

Our constructions are inspired by some Riemannian $n$-disks whose boundaries are ``hard to contract'' through surfaces of controlled $(n-1)$-volume. These Riemannian disks were constructed to answer questions of Gromov \cite{gromov1992asymptotic} and P.~Papasoglu \cite{Papasoglu_ContractingThinDisks} which were special cases of the following general question: Given a Riemannian $n$-disk $D$, can $\partial D$ be homotoped to a point through $D$ while passing through only surfaces whose $(n-1)$-volumes are bounded by a given combination of the geometric parameters of $D$ and $\partial D$?

Riemannian 2-disks whose boundaries are ``hard to contract'' through curves of length bounded in terms of the disk's diameter were constructed by S.~Frankel and M.~Katz. More precisely, for each each $C > 0$ they constructed a Riemannian 2-disk $D^2_C$ such that $\diam(D^2_C), \length(\partial D^2_C) \leq 1$, but every nullhomotopy of $\partial D^2_C$ must pass through a curve longer than $C$ \cite{frankel1993morse}.\footnote{Nevertheless, the boundary of a Riemannian 2-disk $D$ can always be contracted through curves of length at most $\length(\partial D) + 200\diam(D) \max\Big\{1, \ln \frac{\sqrt{\area(D)}}{\diam(D)}\Big\}$, as proven by Liokumovich, Nabutovsky, and Rotman \cite{liokumovich2015contracting}.} Similarly in dimension 3, P.~Glynn-Adey and Z.~Zhu constructed a family of Riemannian 3-disks $D^3_C$, ranging over each $C > 0$, such that $\diam(D^3_C), \vol(D^3_C) \leq 10$ and $\area(\partial D^3_C) = 4\pi$, but every null-homotopy of $\partial D^3_C$ must pass through some surface of area greater than $C$ \cite{GlynnAdeyZhu}. Their construction of $D^3_C$ was inspired by a previous construction by D.~Burago and S.~Ivanov of Riemannian 3-tori with small asymptotic isoperimetric constants \cite{BuragoIvanov_IsoperiConst}.

Our constructions in \cref{thm:theorem} are essentially a combination of $D^3_C$ and a modification of $D^2_C$ by Liokumovich \cite{liokumovich2014surfaces}. Let us compare our constructions with $D^3_C$ and their predecessor constructions by Burago and Ivanov. Burago and Ivanov constructed Riemannian 3-tori that contained 3 intertwining solid tori, each homeomorphic to $D^2 \times S^1$. The metric on each solid torus is a product metric of a short metric on the $S^1$ factor with a Euclidean metric on the $D^2$ factor scaled by a large positive number. In contrast, $D^3_C$ is a Riemannian 3-disk containing two disjoint and linked solid tori $T_1 \cup T_2$, each $T_i$ being homeomorphic to $D^2 \times S^1$. The metric on each $T_i$ is the product of a short metric on the $S^1$ factor and a hyperbolic metric $g$ on the $D^2$ factor with negative curvature of a large magnitude, so that $(D^2,g)$ has small diameter but large area. The metric $M$ we constructed is similar to $D^3_C$, except that in each $T_i$ we replace the hyperbolic metric on the $D^2$ factor with another metric with small diameter and large area.

Our choice of the metric on the $D^2$ factor in each $T_i$ can be motivated from the fact that we are considering maps $F : S^3 \to M$ of nonzero degree and studying the lengths of the loops $F(S^1_{st})$. Each solid torus $T_i$ only ``sees'' part of each loop, namely $F(S^1_{st}) \cap T_i$, which, assuming that the intersection is transverse, is either a closed curve or a union of arcs with endpoints on $\partial T_i$. In other words, $F(S^1_{st}) \cap T_i$ is a relative 1-cycle in $(T_i, \partial T_i)$. These relative 1-cycles form a 2-parameter family ranging over $s$ and $t$; we will demonstrate that some 1-parameter subfamily of those relative 1-cycles will project onto $D^2$ to produce a ``sweepout of $D^2$ by relative 1-cycles,'' a notion that will be defined in \cref{sec:Intro_SweepoutsByCycles}. Proving this involves certain technical complications that we will resolve.

In light of the above, we will construct $M$ so that $D^2$ factor in each $T_i$ has a Riemannian metric that is difficult to sweep out by short relative 1-cycles. This would imply that one of the relative 1-cycles in some $T_i$, and therefore one of the $F(S^1_{st})$'s, has to be long. This metric on $D^2$ will be adapted from Riemannian 2-spheres that are hard to sweep out by short 1-cycles, which were constructed by Liokumovich \cite{liokumovich2014surfaces} following inspiration from \cite{frankel1993morse,liokumovich2013spheres}.

\subsection{Min-max theory for sweepouts by cycles}
\label{sec:Intro_SweepoutsByCycles}

Roughly speaking, if we take the min-max theory for sweepouts by free loops and replace free loops by relative cycles, we obtain a new min-max theory called \emph{Almgren-Pitts min-max theory} that has been central to the study of minimal submanifolds.

More precisely, for any Riemannian $n$-manifold $M$ with boundary, consider the space of \emph{relative flat $k$-cycles in a $M$ with coefficients in $G$}, denoted by $\Cyc_k(M,\partial M; G)$. Intuitively one can think of it as the group of relative singular $k$-cycles in $M$ endowed with a topology where two relative cycles are ``close'' when their difference can be filled by a $(k+1)$-chain of small volume in $M$; formal definitions are available in \cite{FedererFleming_NormalIntegralCurrents,Fleming_FlatChains}. One can define a \emph{sweepout of $M$ by relative $k$-cycles} to be a continuous map $f : N \to \Cyc_k(M, \partial M; G)$ from a simplicial complex $N$ so that $f^*(\iota) \neq 0$, where $\iota \in H^{n-k}(\Cyc_k(M, \partial M; G); G)$ is the \emph{fundamental cohomology class} of $\Cyc_k(M, \partial M; G)$. More generally, for any integer $p \geq 1$, $f$ is called a \emph{$p$-sweepout by cycles} if $f^*(\iota^p) \neq 0$, where $\iota^p$ denotes the $p^\text{th}$ cup power. A min-max value called the \emph{$p$-width} can be defined as:
\begin{equation}
    \width_p^k(M; G) = \inf_{\substack{
        f : N \to \Cyc_k(M, \partial M; G) \\
        f^*(\iota^p) \neq 0
    }} \sup (\vol_k \circ f).
\end{equation}
Henceforth we will write $\width_p^k(M)$ to denote $\width_p^k(M; \Z)$ when $M$ is orientable, and $\width_p^k(M; \Z_2)$ otherwise. The $p$-widths form a non-decreasing sequence: $\width_1^k(M) \leq \width_2^k(M) \leq \width_3^k(M) \leq \dotsb$

The $p$-widths of a closed Riemannian manifold $M$ are realized as the volumes of minimal submanifolds in $M$, which may contain a ``small'' singular set. The study of $p$-widths have led to existence proofs for minimal submanifolds and the solutions of several conjectures about minimal submanifolds. Almgren-Pitts Min-max theory and some its applications to these conjectures are surveyed in \cite{Marques_MinimalSurfacesSurvey,Neves_MinMaxApplications}.

Since every free loop in $M$ is an integral 1-cycle, we have $\width_1^1(M) \leq \lambda(M)$. Like $\lambda(M)$, curvature-free bounds on the $p$-widths of $M$ in terms of geometric parameters of $M$ are much better understood when $\dim M = 2$. For every closed Riemannian surface $S$ and $p \geq 1$ we have $\scg(S) \leq \width_p^1(S)$ due to the recent work of O.~Chodosh and C.~Mantoulidis \cite{ChodoshMantoulidis_pWidthsGeodesics}. F.~Balacheff and Sabourau proved that $\width_1^1(S) \leq 10^8\sqrt{\genus(S)+1}\sqrt{\area(S)}$ \cite{BalacheffSabourau_Diastoles}. However, Liokumovich proved that $\width_1^1(S)$ (and, by extension, $\width_p^1(S)$ for all $p$) cannot be bounded solely in terms of $\diam(S)$, answering a question of Sabourau \cite{Sabourau_ShortestClosedGeodesicS2_Diameter}: the counterexamples from \cite{liokumovich2013spheres} can be adapted into a family of Riemannian surfaces $S$ for which $\diam(S) = 1$ but whose values of $\width_1^1(S)$ grow without bound \cite{liokumovich2014surfaces}. This exemplifies the potential for bounds on $\lambda$, and counterexamples to conjectured bounds, to shed light on the $p$-widths, and by extension the geometry of minimal submanifolds.

We adapted Liokumovich's Riemannian surfaces of small diameter but large 1-width into Riemannian 2-disks with small diameter but large 1-width. These Riemannian 2-disks are hard to sweep out by short relative 1-cycles, and they serve as ingredients in our construction, as explained at the end of \cref{sec:MotivatingOurConstruction}.

For manifolds of dimension at least 3, bounds on some of their $p$-widths have been established under certain curvature assumptions \cite{Gromov_FillingRiemannianManifolds,Guth_WidthVolIneq,MarquesNeves_Widths_PosRicci,ParkerLiokumovich_Width_PosRicci,Sabourau_Waist_PosRicci,LiokumovichMaximo_Waist_PosScalar,LiokumovichZhou_Waist_PosRicci}. However, there are no known curvature-free bounds on the $p$-widths of manifolds of dimension 3 and above in terms of geometric parameters such as diameter and volume. In fact, for closed 3-manifolds $M$ and any $p$, $\width_p^2(M; \Z_2)$ cannot be bounded in terms of $\diam(M)$ and $\vol(M)$. This can be shown in two steps: roughly speaking, $M$ can be ``cut'' into two regions of equal volume by a hypersurface of area $\width_1^2(M; \Z_2)$.\footnote{This follows due to an argument adapted from \cite[p.~396]{liokumovich2014surfaces}.} On the other hand, Papasoglu and E.~Swenson constructed, for each $C > 0$, a Riemannian 3-sphere of diameter and volume at most 1 for which any such ``cutting hypersurface'' must have area greater than $C$ \cite{PapasogluSwenson_Expander}.\footnote{Riemannian 3-disks that are ``hard to cut'' in this manner were independently constructed by Glynn-Adey and Zhu \cite{GlynnAdeyZhu}, except that they only studied cutting hypersurfaces that were embedded disks.}

The preceding argument would not apply to widths $\width_p^k(M)$ where $k \leq \dim M - 2$. Nevertheless, our \cref{thm:theorem} may serve as a first step towards proving that the $p$-widths $\width_p^1$ of Riemannian 3-spheres cannot be bounded in terms of diameter and volume.

\subsection{Min-max theory for sweepouts by slicing}

Our result may also offer insights in the study of \emph{waists}, which are another class of min-max values that were defined by Gromov in \cite{Gromov_FillingRiemannianManifolds}. They arise from sweepouts of Riemannian $n$-manifolds with boundary $M$ by relative $k$-cycles that are obtained as ``slices'' of $M$. More precisely, given a continuous map $M \to \R^{n-k}$, its fibers can be considered as ``slices'' of $M$. We can define the \emph{$k$-waist} of orientable $n$-manifolds $M$ to be
\begin{equation*}
    \waist_k(M) = \inf_{\sigma : M \to \R^{n-k}} \sup_{t \in \R^{n-k}} \vol_k(\sigma^{-1}(t)),
\end{equation*}
where the infimum is taken over maps $\sigma$ whose fibers are Lipschitz $k$-cycles, such that the map $\R^{n-k} \to \Cyc_k(M, \partial M; \Z)$ given by $t \mapsto \sigma^{-1}(t)$ is continuous.

The known bounds on waists in terms of the geometric parameters of $M$ follow a pattern similar to the other min-max values that we introduced earlier. By the definitions we have $\width_1^k(M) \leq \waist_k(M)$. For a Riemannian 2-sphere $(S^2,g)$, Liokumovich proved that $\waist_1(S^2,g) \leq 52\sqrt{\area(S^2,g)}$ \cite{Liokumovich_Waist_Sphere}.\footnote{A similar bound on $\waist_1$ for other closed Riemannian surfaces is implicit via the \emph{monotone sweepouts} involved in the proof of \cite[Theorem~1.1]{Liokumovich_ImplicitWaists_Surfaces}.} For manifolds of dimension at least 3, some bounds on waists have been obtained under curvature assumptions \cite{LiokumovichZhou_Waist_PosRicci,LiokumovichMaximo_Waist_PosScalar,Sabourau_Waist_PosRicci}.

There have been more results on $\waist_k(M)$ when $k = \dim M - 1$; the aforementioned Riemannian 3-spheres that are ``hard to cut'' by hypersurfaces demonstrate that $\waist_2(M)$, which is bounded from below by $\width_1^2(M)$, cannot be bounded in terms of $\diam(M)$ and $\vol(M)$ when $\dim M = 3$. On the flipside, it is currently an open question whether $\waist_k(M)$ can be bounded in terms of the geometric parameters of $M$ when $\dim M \geq 3$ and $k \leq \dim M - 2$. Within this context, L.~Guth conjectured that when $M$ is a Riemannian 3-torus, then $\waist_1(M) \leq C\sqrt[3]{\vol(M)}$ for some constant $C$ \cite[Section~7]{Guth_MetaphorsSystolic}. Our \cref{thm:theorem} could serve as a first step towards disproving this conjecture.

Our result can be contrasted with a recent result of Nabutovsky, Rotman, and Sabourau which proves bounds on another min-max value related to waists. To define this min-max value, for each closed Riemannian $n$-manifold $M$ and an integer $k \geq 0$, consider a continuous map $F : N \to M$ of nonzero degree from a closed $n$-dimensional pseudomanifold $N$. ``Slice'' $N$ similar to before using a continuous map $\sigma : N \to K$ to a finite simplicial complex $K$ whose fibers $\sigma^{-1}(t)$ are $k$-dimensional simplicial complexes. Then define
\begin{equation*}
    W_k(M) = \inf_{\substack{
        \sigma : N \to K \\
        F : N \to M \\
        \deg f \neq 0
    }} \: \sup_{t \in K} \: \vol_k(F(\sigma^{-1}(t))),
\end{equation*}
where the infimum range over all $\sigma$ and $F$ that satisfy the stipulated criteria.

These min-max values are related to waists via $W_k(M) \leq \waist_k(M)$. The bounds $W_k(M) \leq c_n\sqrt[n]{\vol(M)}$ and $W_k(M) \leq c_n'\diam(M)$ were proven by Nabutovsky, Rotman, and Sabourau for some dimensional constants $c_n$ and $c_n'$ \cite{NabutovskyRotmanSabourau_SweepoutsSlicing}. However, it remains to be seen whether $W_k(M)$ is attained as the volume of some minimal object. Since each $F(\sigma^{-1}(t))$ may not be a loop or a cycle (see \cite[Example~1.2]{NabutovskyRotmanSabourau_SweepoutsSlicing}), the min-max theory for minimal submanifolds does not fit here. For the same reason, there is no direct comparison between $W_1(M)$ and $\lambda(M)$ or $\width_p^1(M)$.

\subsection{Orthogonal geodesic chords and min-max theory for sweepouts by paths}
\label{sec:Intro_SweepoutsByPaths}

The techniques we used to prove \cref{thm:theorem} can be adapted for yet another min-max theory that arises from using sweepouts by paths instead of free loops or 1-cycles. This has led to a result about the geometry of orthogonal geodesic chords, which were defined earlier as geodesics in a Riemannian manifold $M$ that meet a fixed submanifold $N$ orthogonally at its endpoints. Orthogonal geodesic chords are 1-dimensional analogues of \emph{free boundary minimal submanifolds}, which are submanifolds $P \subset M$ with vanishing mean curvature such that $\partial P \subset N$ and $\partial P$ meets $N$ orthogonally.\footnote{Some results about the existence and regularity of free boundary minimal submanifolds are surveyed in \cite{Li_FreeBoundaryMinimalSurfaces}.} Orthogonal geodesic chords are also related to \emph{brake orbits}, special types of periodic orbits in certain Hamiltonian systems \cite{GiamboGiannoniPiccione_OGC}.

When $\dim N = 0$, an orthogonal geodesic chord is simply a geodesic with specified endpoints. The existence of such geodesics and bounds on their length were studied in \cite{Serre_GeodesicLoops,Schwartz_GeodesicArcs,NabutovskyRotman_GeodesicsQuadBound,NabutovskyRotman_GeodesicLoopsLinearBound,Cheng_GeodesicsLinearBound,Beach_SimpleGeodesicLoops}. Lyusternik and L.~Schnirelmann proved that every convex domain $M \subset \R^n$ with boundary $N$ contains $n$ orthogonal geodesic chords \cite{LyusternikSchnirelmann_OrthogonalGeodesicChords}. W.~Bos extended this existence result to Riemannian $n$-disks $M$ with convex boundary $N$ \cite{Bos_OrthogonalGeodesicChords}. For Riemannian 2-disks $M$ with strictly convex boundary $N$, J.~Hass and P.~Scott \cite{HassScott_OrthogonalGeodesicChords} and D.~Ko \cite{Ko_FreeBoundaryGeodesic} showed that one can even arrange the orthogonal geodesic chords to be \emph{simple}, that is, avoid self-intersection.

The existence and geometry of orthogonal geodesic chords may be probed using min-max techniques as follows. Let $\Omega_N M$ denote the space of piecewise smooth paths in $M$ that start and end on $N$, topologized as in \cite[p.~88]{Milnor}. Let $d \geq 1$ be the smallest degree for which $\pi_d(\Omega_N M, \Lambda^0 N) \neq 0$. Then define
\begin{equation}
    \lambda_\rel(M, N) = \inf_{\substack{
        f : (D^d, \partial D^d) \to (\Omega_N M, \Lambda^0 N) \\
        0 \neq [f] \in \pi_d(\Omega_N M, \Lambda^0 N)
    }} \sup (E \circ f),
\end{equation}
where $E$ is the \emph{energy} of a path, $E(\alpha) = \int_0^1 \norm{\alpha'(t)}^2\,dt$. (Using length instead of energy would not significantly affect the resulting min-max theory, because the H\"older inequality relates length to energy.) X.~Zhou proved that when $M$ is a complete and homogeneously regular Riemannian manifold and $N$ is a closed submanifold such that $\pi_1(M, N) = 0$ and $\pi_2(M,N) \neq 0$, then $\lambda_\rel(M, N) > 0$ is the energy of an orthogonal geodesic chord \cite{Zhou_FreeBoundaryGeodesics}.

Recent results have estimated the lengths of orthogonal geodesics chords in various spaces $M$ and $N$, where $N$ may or may not be $\partial M$. When $M$ is a Riemannian 2-disk with convex boundary $N$, I.~Beach proved that $M$ contains at least two distinct simple orthogonal geodesic chords whose lengths are bounded by $f(\diam(M), \area(M), \length(\partial M))$ for some function $f$ \cite{Beach_SimpleOGC}. In addition, recent work by Beach, H.\,C.~Peruyero, E.~Griffin,
M.~Kerr, Rotman, and C.~Searle implies that when $M$ is a closed Riemannian manifold and $N$ is an analytic 2-sphere embedded in $M$, then $M$ contains an orthogonal geodesic chord whose length is bounded by $f'(\dim M, \diam(M), \diam(N), \area(N))$ for some function $f'$ \cite{BeachEtAl_OGC}. Both of these results were obtained by proving that either $\lambda_\rel(M,N)$ is bounded by the relevant function ($f$ or $f'$) due to the existence of a sweepout of $M$ by curves of energy (equivalently, length) bounded by that function, or else the obstruction to the existence of such a sweepout is an orthogonal geodesic chord of length bounded by that function.

When $M$ is a Riemannian 3-sphere and $N$ is an embedded sphere, we proved that $\lambda_\rel(M,N)$ cannot be bounded by any function of $\vol(M)$, $\diam(M)$, and $\diam(N)$.
\begin{theorem}
    \label{thm:OrthoGeodesicChords}
    For any $E > 0$, there exists a Riemannian 3-sphere $M$ of diameter and volume at most 1 that contains an embedded 2-sphere $N$ of diameter at most 1 such that $\lambda_\rel(M, N)$. $M$ also contains an embedded circle $\gamma$ of length at most 1 such that $\lambda_\rel(M, \gamma) > E$.
\end{theorem}

\section{The Construction of our 3-Spheres}
\label{sec:Construction}

\begin{figure}[h]
    \centering
    \begin{tabular}{cc}
        (a) & (b)
        \\
        \includegraphics[width=0.3\linewidth]{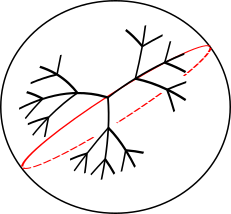} & \includegraphics[width=0.3\linewidth]{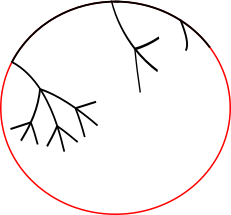}
        \\
        (c) & (d)
        \\
        \includegraphics[width=0.3\linewidth]{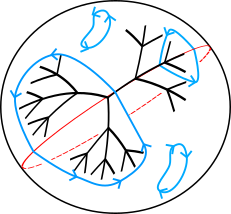} & \includegraphics[width=0.3\linewidth]{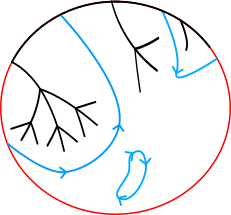}
    \end{tabular}
    \caption{(a) A 2-sphere $S_h$ of small diameter but large width from 
    \cite{liokumovich2014surfaces}. (b) A Riemannian disk $D_h$ whose double is $S_h$. (c) A 1-cycle (blue) in $S_h$, which is the double of the relative 1-cycle $A$  (blue) in $D_h$, shown in (d).}
    \label{fig:HalfSphereRelativeCycle}
\end{figure}

For any $L > 0$, Liokumovich constructed a Riemannian 2-sphere $S_h$ of diameter at most 1 for which $\width_1^1(S_h) > L$ \cite[pg.~2]{liokumovich2014surfaces}. $S_h$ is roughly constructed as follows. Take a unit disk $B$ in a hyperbolic plane, and embed a regular ternary tree with unit edge length and height $h$ in $\R^2$. Glue $B$ to a small tubular neighbourhood of the tree in $\R^2$ by identifying their boundaries. (The curvature of $B$ must be chosen such that the boundaries have the same length.) The result is $S_h$ (see \cref{fig:HalfSphereRelativeCycle}(a)). 

The metric of $S_h$ has the symmetries of an equilateral triangle, generated by reflections and rotations. A plane of reflection cuts $S_h$ into two isometric Riemannian disks, one of which we denote by $D_h$ (see \cref{fig:HalfSphereRelativeCycle}(b)). It can be verified that $\diam D_h \leq \frac12$.

The key property of $D_h$ is that it has a small diameter but large width:

\begin{lemma}
    \label{lem:DiskLargeWidth}
    For any $C > 0$, there exists some $h$ such that $\width_1^1(D_h) > C$. 
\end{lemma}
\begin{proof}
    Consider a sweepout by relative 1-cycles $f : N \to \Cyc_1(D_h, \partial D_h; \Z)$; unpacking the definition, this implies that for some loop $g : S^1 \to N$, the gluing homomorphism sends $g \circ f$ to $m$ times of the fundamental class of $D_h$, for some $m \neq 0$. Consider the map $f' : N \to \Cyc_1(S_h; \Z)$ where the 1-cycle $f'(x)$ (see \cref{fig:HalfSphereRelativeCycle}(c)) is the double of the relative 1-cycle $f(x)$ (see \cref{fig:HalfSphereRelativeCycle}(d)), which is obtained by subtracting a reflected copy of $f(x)$ from $f(x)$. (The subtraction ensures that the orientations match up where the two copies of $f(x)$ meet.) Then it can be verified that the gluing homomorphism sends $f' \circ g$ to $\pm m$ times of the fundamental class of $S_h$. The lemma then follows from the fact that $S_h$ can have arbitrarily large $\width_1^1$.
\end{proof}

For $i=1,2$, define $T_i = D_h \times S^1$ to be solid tori, and let $g'_h$ be the product metric of $D_h$ with a sufficiently short metric on $S^1$. Embed the disjoint union $T_1 \cup T_2$ into $S^3$ in the manner of a Hopf link (see \cref{fig:linking_core_curve}), and denote their union by $LT$. Cover $S^3$ by two open sets $U_1$ and $U_2$ so that $LT \subset U_1 \setminus U_2$, and let $\phi_1, \phi_2$ be a partition of unity subordinate to this cover. Extend $g'_h$ over $U_1$ via a smooth bump function. Our Riemannian $3$-sphere is then $M_h = (S^3,g_h)$, where $g_h = g'_h\phi_1 + g_0\phi_2$ and $g_0$ is the metric of a round sphere of sufficiently small radius. It can be verified that $\vol(M_h), \diam(M_h) \leq 1$.

\section{Plan of the Proof of \cref{thm:theorem}} \label{Plan}

In order to discuss the plan, we first make some notation clear. Let $\pi_i: T_i \ra D_h$ be a projection onto the first factor. Recall the definition of $S^1_{st}$ from the beginning of the Introduction. Define $S^2_t = S^3 \cap (\R^3 \times \{t\})$.

Consider a map $F : S^3 \to M_h$ of nonzero degree. After perturbing $F$, each $X_i = F^{-1}(\partial T_i)$ will be a smooth manifold containing a family of 1-cycles $X_i \cap S^2_t$. Their images $F(X_i \cap S^2_t)$ will form a sweepout of $\partial T_i$. An analysis of $\deg F|_{X_i}$ and the following lemma of algebraic topology guarantees the existence of a circle $C$ in some $X_i \cap S^2_t$ that maps to a non-contractible loop $F(C)$ on $\partial T_i$:

\begin{lemma}
\label{lem:NonContractibleOnMorseFunctionPreimage}
    Let $X$ be a closed, orientable, and connected surface of genus at least 1 and consider a degree nonzero map $\alpha : X \to S^1 \times S^1$. Then for any Morse function $\phi : X \to \R$, there exists some $t \in \R$ and a circle $C$ embedded in $\phi^{-1}(t)$ such that $\alpha|_C$ is not nullhomotopic.
\end{lemma}

This result follows from \cite[Lemma~2.2]{GlynnAdeyZhu}. Nevertheless, at the end of this section we will sketch a proof for the case where $\alpha$ is a diffeomorphism, in order to articulate the fundamental reason why it is true.

\begin{figure}[h]
    \centering
    \includegraphics[width=0.5\linewidth]{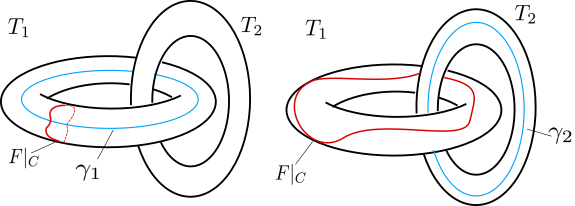}
    \caption{When $F|_C : C \to \partial T_1$ (red curve) is not nullhomotopic, it has a nonzero linking number with either the core curve $\gamma_1$ of $T_1$ (left picture) or the core curve $\gamma_2$ of $T_2$ (right picture).}
    \label{fig:linking_core_curve}
\end{figure}

As shown in \cref{fig:linking_core_curve}, we find that $F|_C$ has a nonzero linking number with the core curve of either $T_1$ or $T_2$. As a consequence of \cref{lem:winding}, there will exist a surface $\Sigma \subset S^2_t$ such that $F(\Sigma)$ is contained within some $T_i$. In addition, $\Sigma$ will be partitioned into the curves $\Sigma \cap S^1_{st}$, and the projections of the curves $F(\Sigma \cap S^1_{st})$ onto $D_h$ via $\pi_i$ will give a sweepout of $D_h$ by relative 1-cycles. \cref{fig:EssentialDisk} illustrates $\Sigma$ and $F(\Sigma)$, where $\Sigma$ corresponds to $\Sigma_1$ in the figure.

One of the curves in this sweepout by relative 1-cycles must be long, by \cref{lem:DiskLargeWidth}. As these relative 1-cycles are orthogonal projections of the curves $F(\Sigma \cap S^1_{st})$, it must be the case that some $F(S^1_{st})$ is long as well.

\begin{proof}[Proof sketch of \cref{lem:NonContractibleOnMorseFunctionPreimage} when $\alpha$ is a diffeomorphism]
    Since $\alpha$ is a diffeomorphism, let us identify $X$ with $S^1 \times S^1$. Let the critical values of $\phi$ be $t_0 < t_1 < \dotsb < t_m$. For each $i$, $\phi^{-1}(t_i)$ is a disjoint union of wedge sums of circles. Denote those circle wedge summands by $C_{i,1}, \dotsc, C_{i,k_i}$ (see \cref{fig:SphericalDecomposition}(a)). Let us assume that each circle $C_{i,j}$ is contractible in $X$ and derive a contradiction. In other words, we assume that each $C_{i,j}$ is the boundary of some 2-chain $D_{i,j}$ in $X$ that is the continuous image of a 2-disk (see \cref{fig:SphericalDecomposition}(c)). 
    
    We can decompose the fundamental cycle of $X$ into the sum of 2-chains $A_i = \phi^{-1}([t_{i-1}, t_i])$ for $i = 1,\dotsc, m$ (see \cref{fig:SphericalDecomposition}(b)). Each $\partial A_i$ is composed of circles $C_{i',j'}$ (for $i' = i-1, i$), and for each of these circles we can glue in $D_{i',j'}$. After ``capping off'' each circle in $\partial A_i$, $A_i$ becomes a 2-cycle that is the union of images of 2-spheres. To illustrate, in \cref{fig:SphericalDecomposition}(b), we can see that $A_1$ and $A_3$ have been capped off into spheres, while $A_2$ has been capped off into the union of two spheres.
    
    \begin{figure}[h!]
        \centering
        \includegraphics[scale=0.15]{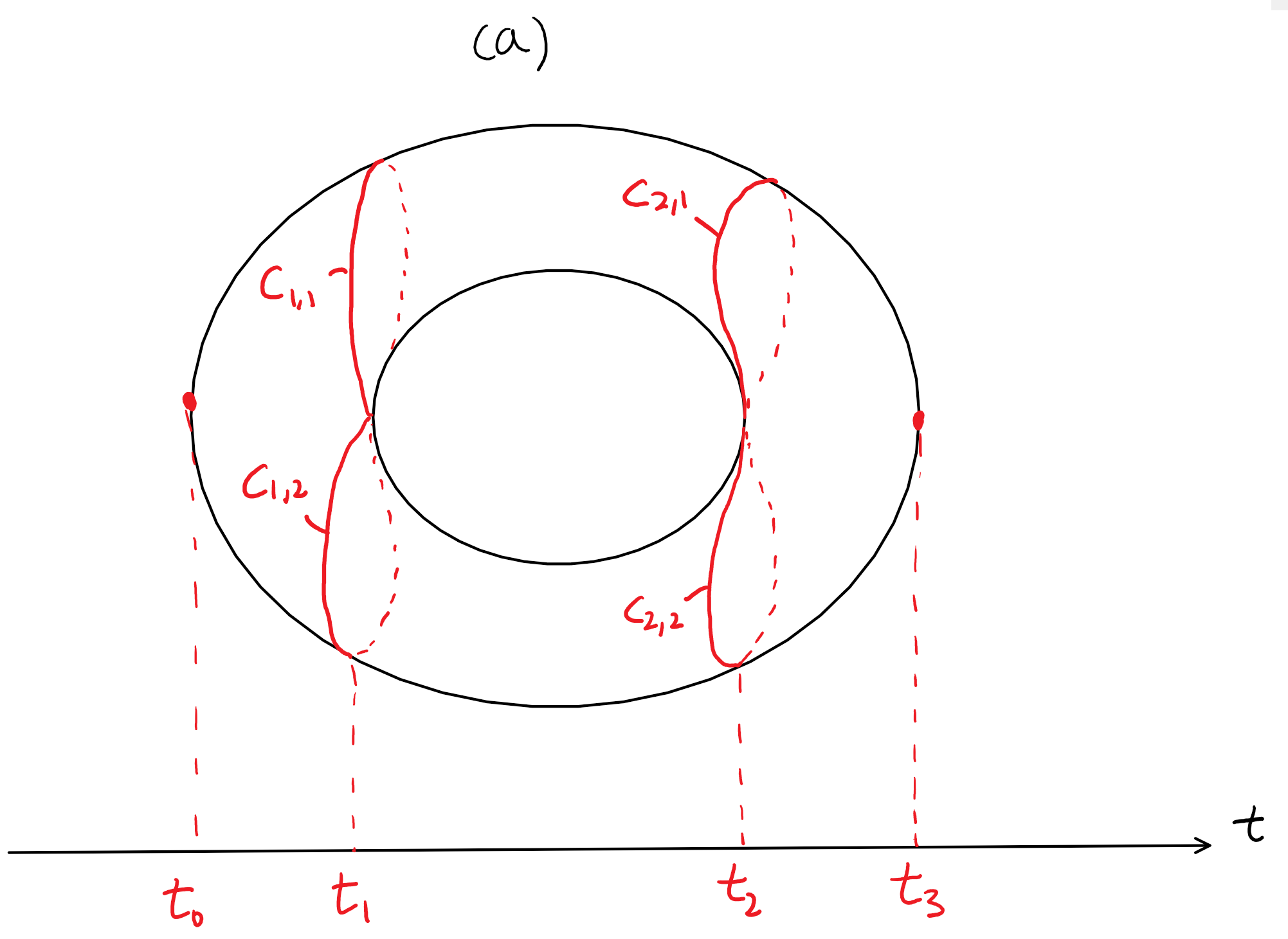}
        \\[2ex]
        \includegraphics[scale=0.15]{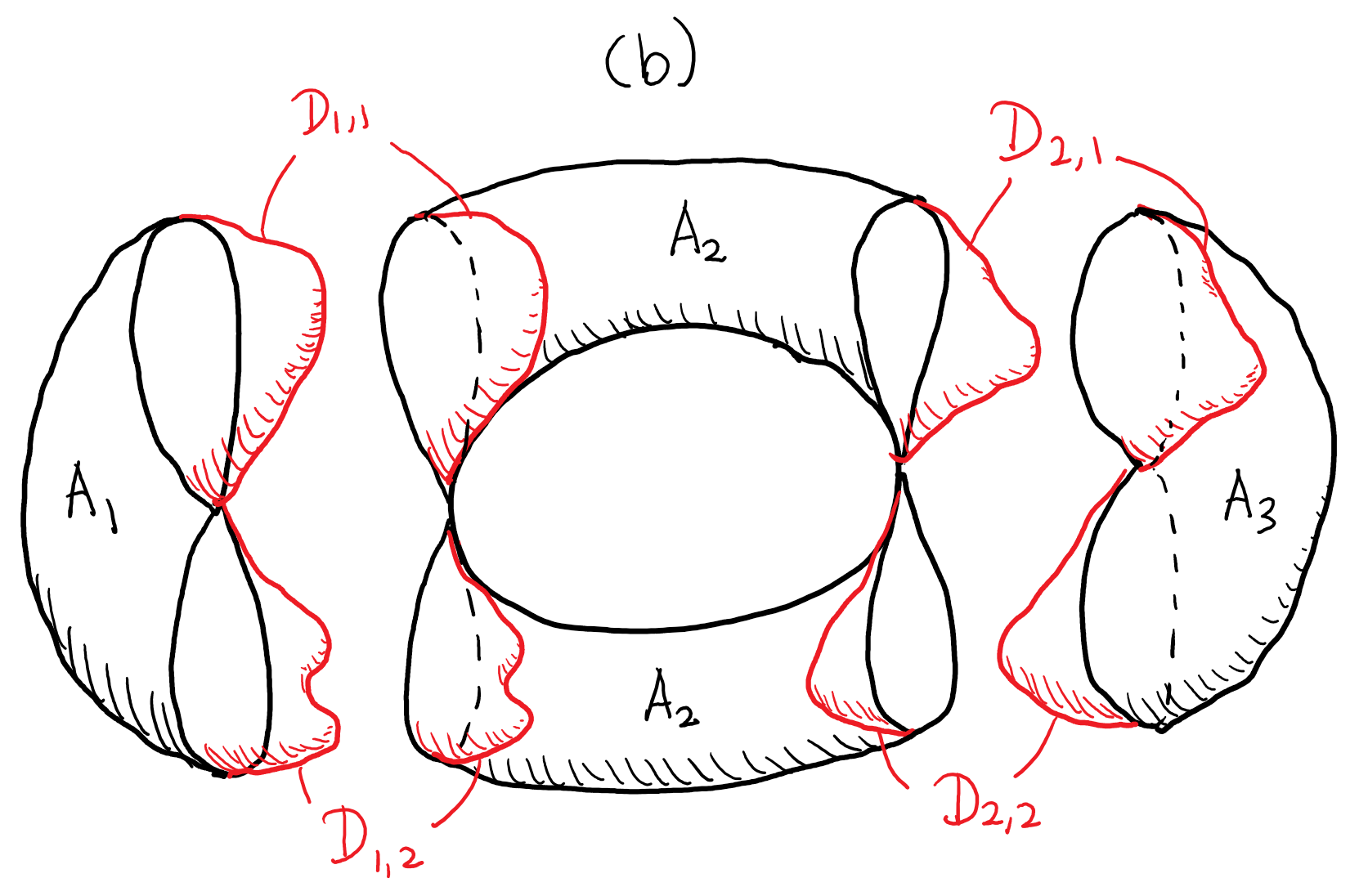}
        \quad
        \includegraphics[scale=0.15]{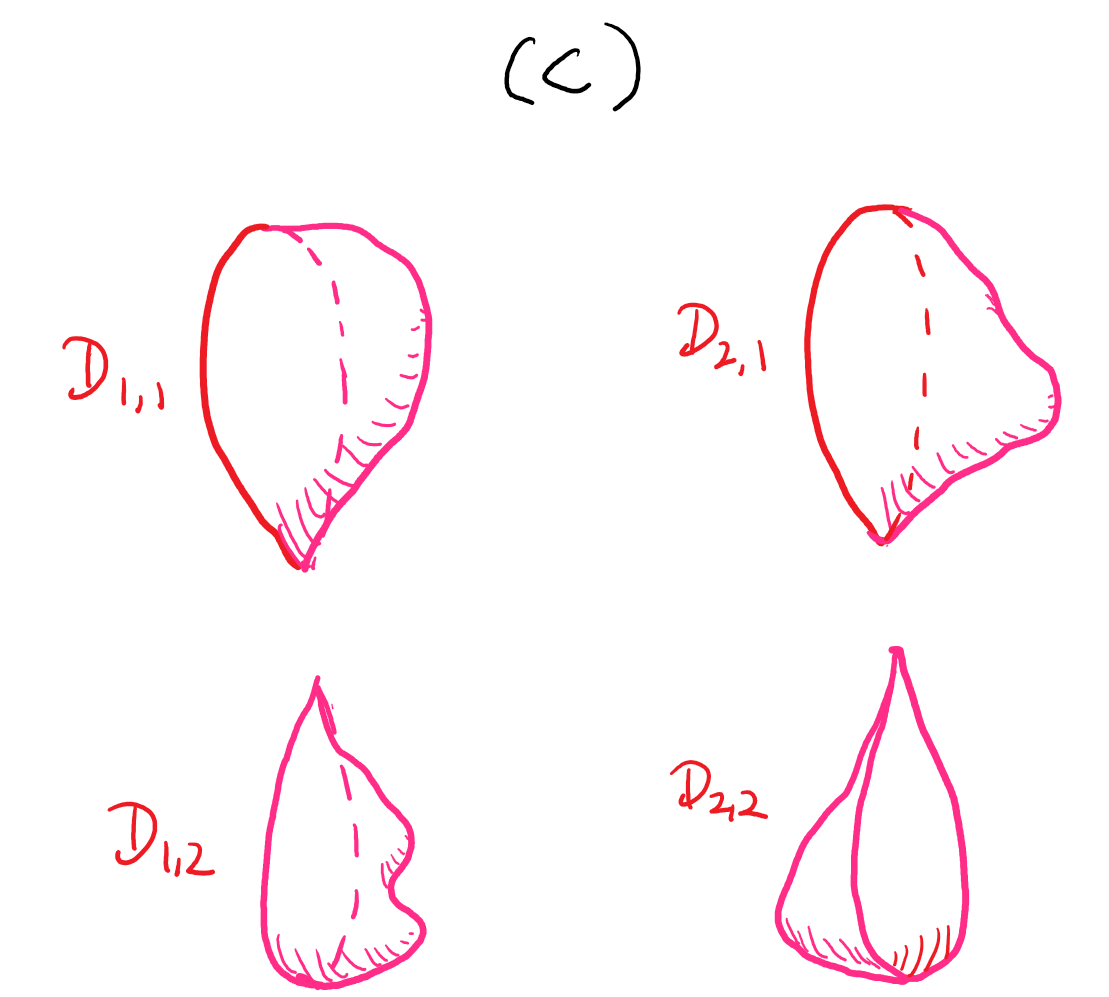}
        \caption{(a) Level sets of a Morse function on a 2-torus $X$. (b) ``Capping off'' each circle in each level set produces a decomposition of the fundamental cycle of $X$ into a sum of images of spheres. (c) The images of disks used for ``capping off.''}
        \label{fig:SphericalDecomposition}
    \end{figure}
    
    However, as $\pi_2(X) = 0$, all of these images of 2-spheres are null-homologous. If we add up all of these images of spheres, then the caps $D_{i,j}$ cancel each other out and the result is the fundamental cycle of $X$. Thus we have expressed the fundamental cycle as a sum of null-homologous 2-cycles, which gives a contradiction.
\end{proof}

\section{Main Result}
\label{sec:MainResult}

We begin by proving that continuous maps $S^3 \to M$ of nonzero degree can be perturbed to ``geometrically nice'' maps.

\begin{lemma}
    \label{lem:morse}
    Consider any map $F : S^3 \ra M$ of nonzero degree such that every curve $F(S^1_{st})$ is shorter than $L$ for some $L > 0$. Then for any $\delta > 0$ and any closed submanifold $Y \subset M$, $F$ is homotopic to a smooth map $\hat{F}$ transverse to $Y$ such that:
    \begin{enumerate}
        \item For any $s$ and $t$, the lengths of $F(S^1_{st})$ and $\hat{F}(S^1_{st})$ differ by at most $\delta$.

        \item The sets $\hat{F}^{-1}(Y) \cap S^2_t$ are the level sets of some Morse function on $\hat{F}^{-1}(Y)$.
    \end{enumerate}
\end{lemma}
\begin{proof}
    We can approximate $F$ by a sequence of smooth maps $F_i : S^3 \to M$ for $i = 1, 2, \dotsc$ that are transverse to $Y$, so that each $X_i = F_i^{-1}(Y)$ is a smooth manifold. It can be verified from the standard arguments for such approximations (e.g. in \cite{Lee_IntroSmoothManifolds}), together with the fact that every $F(S^1_{st})$ is shorter than $L$, that each $F_i$ can be chosen to satisfy (1). Thinking of $X_i$ as a submanifold of $\R^4$, the function $d_p : X_i \to \R^4$ defined by $d_p(x) = \norm{x - p}^2$ is Morse for generic $p \in \R^4$ \cite[p.~36]{Milnor}. Let us pick $p \approx (0,0,0,2)$; then the rotational symmetry of $S^3$ implies that the level sets of $d_p$ are the intersections of $X_i$ with hyperplanes orthogonal to $p$. Therefore we may modify $F_i$ and $d_p$ by precomposing them with an isometry $S^3 \to S^3$ that is close to the identity, until those level sets become $X_i \cap S^2_t$. Finally, we may choose $\hat{F}$ to be $F_i$ for sufficiently large $i$.
\end{proof}

Thus we may replace $F$ by its perturbation $\hat{F}$. Henceforth we will assume that $F$ satisfies the properties in \cref{lem:morse}. Let $X = F^{-1}(\partial T_1 \cup \partial T_2)$. Then like any level set of a Morse function on a surface, each $X \cap S^2_t$ is a disjoint union of circles and at most one figure eight (a wedge sum of two circles).

The proofs of the following two lemmas were inspired by the proof of \cite[Lemma~2.3]{GlynnAdeyZhu}. However, that proof contained an inaccuracy, so we give our own proofs.

\begin{figure}
    \centering
    \includegraphics[width=12cm]{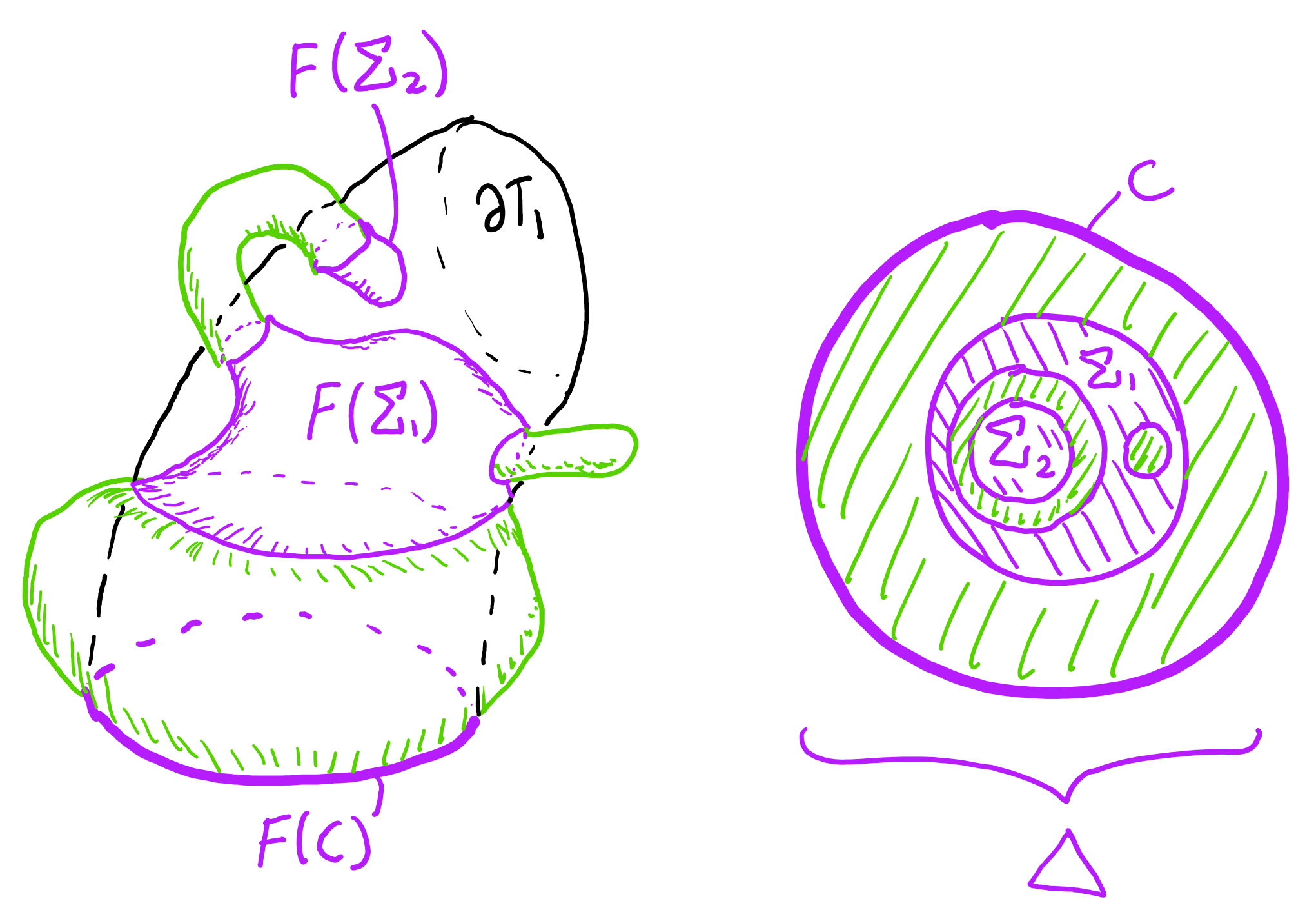}
    \caption{An illustration of some elements from the proof of \cref{lem:SurfaceExistence}.}
    \label{fig:EssentialDisk}
\end{figure}

\begin{lemma}
    \label{lem:SurfaceExistence}
    Suppose that one of the circles $C$  embedded in $X \cap S^2_t$ is such that $F|_C$ has nonzero linking number $m$ with the core curve of $T_i$ for some $i = 1,2$. Then $F^{-1}(T_i) \cap S^2_t$ contains a surface $\Sigma$ such that $(\pi_i \circ F)_*:H_2(\Sigma, \partial \Sigma) \ra H_2(D_h, \partial D_h)$ is a nonzero map.
\end{lemma}
\begin{proof}
    Without loss of generality, assume $F|_C$ has linking number $m \neq 0$ with the core curve of $T_1$. $C$ then bounds a disk $\Delta$ in $S^2_t$ such that $F(\Delta)$ has intersection number $m$ with the core curve of $T_1$.\footnote{With reference to \cref{rem:Generalization} about generalizations of \cref{thm:theorem}, this is the part of our proof of \cref{thm:theorem} that requires the Jordan curve theorem.} $\Delta \cap F^{-1}(T_1)$ is a union of disjoint surfaces $\Sigma_1 \cup \dotsb \cup \Sigma_k$, each corresponding to an intersection number $m_i$ between $F(\Sigma_i)$ and the core curve of $T_1$. This is reflected in \cref{fig:EssentialDisk}. Then $m = m_1 + \dotsb + m_k \implies m_i \neq 0$ for some $m_i$. Thus $\Sigma_i$ is the desired surface.
\end{proof}

We then use \cref{lem:SurfaceExistence} to prove the following statement:
 
\begin{lemma}
\label{lem:winding}
For some $t$ and for some $i = 1,2$, $F^{-1}(T_i) \cap S^2_t$ contains a surface $\Sigma$ such that $(\pi_i \circ F)_*:H_2(\Sigma, \partial \Sigma) \ra H_2(D_h, \partial D_h)$ is a nonzero map.
\end{lemma}

\begin{proof}
Let $X_1 = F^{-1}(\partial T_1)$. By \cref{lem:morse} we can assume that $X_1 \cap S^2_t = \phi^{-1}(t)$ for some Morse function $\phi : X_1 \to \R$. After choosing an orientation on $\partial T_1$ and equipping $X_1$ with the preimage orientation, it can be verified that the map $F|_{X_1} : X_1 \to \partial T_1$ has degree equal to $\deg F \neq 0$. (To see this, note that the preimages of a regular value of $F$ have neighbourhoods that are foliated by the $S^2_t$'s.) Consequently, for some connected component $X_1'$ of $X_1$, $\deg F|_{X_1'} \neq 0$. Since $\pi_2(\partial T_2) = 0$, the genus of $X_1'$ must be at least 1. Thus \Cref{lem:NonContractibleOnMorseFunctionPreimage} gives some $t \in \R$ and a circle $C$ embedded in $\phi^{-1}(t)$ such that $F|_C$ is not contractible in $\partial T_1$.

In other words, $F_*: \pi_1(C) \iso \Z \ra \pi_1(\partial T_1)$ is a nonzero map with $F_*(1) = (n_1,n_2)$, where the first factor is a multiple of the generator homotopic to $\partial D_h$ and the second factor is a multiple of the generator homotopic through $T_1$ to its core curve.

Then if $n_1 \neq 0$, we have that  $F|_C$ has nonzero linking number with the core curve of $T_1$, and so we may apply \cref{lem:SurfaceExistence} to obtain $\Sigma$ as needed.

If $n_1 = 0$, then $n_2 \neq 0$ where $F|_C$ has a nonzero linking number with the core curve of $T_2$, and so we again apply \cref{lem:SurfaceExistence} to obtain $\Sigma$.
\end{proof}

We have now proven the existence of a surface $\Sigma$ on $S^2_t$ that maps into one of the tori in a ``nice'' way. We use this property to define a continuous family of 1-cycles on $D_h$ and show that one of those 1-cycles must be long. This will eventually imply that one of the $F(S^1_{st})$'s must also be long, leading to a proof of \cref{thm:theorem}.

\begin{proof}[Proof of \cref{thm:theorem}]
Consider any $L > 0$. Consider the Riemannian 3-sphere $M_h = (S^3, g_h)$ that was constructed in \cref{sec:Construction}, with $h$ chosen such that $\width_1^1(D_h) > L + 1$, as in \cref{lem:DiskLargeWidth}. Suppose for the sake of contradiction that for some continuous map $F : S^3 \to M_h$, every curve $F(S^1_{st})$ is shorter than $L$. Recall that by applying a perturbation, we may assume that $F$ satisfies the properties of \cref{lem:morse}, for some $\delta < \frac12$. We will arrive at a contradiction by proving that some $F(S^1_{st})$ (for the perturbed $F$) is longer than $L + 1$.

Applying \cref{lem:winding}, we obtain a surface $\Sigma$ in $S^2_t$ for some $t$ such that $F(\Sigma) \subset T_i$ for some $i$, and $(\pi_i \circ F)_* : H_2(\Sigma, \partial\Sigma) \to H_2(D_h, \partial D_h)$ is a nonzero map. With this, we may define a continuous family of relative 1-cycles $K :[-1,1] \ra \Cyc_1(D_h, \partial D_h; \Z)$ by $K(s) = \pi_i(F(S^1_{st} \cap \Sigma))$. Since $K(\pm1) = 0$, we can think of $K$ as a map $S^1 \to Z_1(D_h, \partial D_h)$.

We will show that $K$ gives a sweepout of $D_h$ by relative 1-cycles by appealing to the Almgren isomorphism theorem, which implies that there is a natural ismorphism $\Gamma : \pi_1(\Cyc_1(D_h, \partial D_h; \Z)) \xrightarrow{\iso} H_2(D_h, \partial D_h)$ \cite{Almgren_AlmgrenIso}.\footnote{A modern proof of the Almgren isomorphism theorem is available in \cite{GuthLiokumovich_ParamIneq}.} $\Gamma$ is induced by ``gluing'' a family of relative 1-cycles in $D_h$ into a relative 2-cycle in $D_h$ (see \cref{fig:GluingHomomorphism}). Gluing together the relative 1-cycles $\pi_i(F(S^1_{st} \cap \Sigma))$ gives $\pi_i(F(\Sigma))$, which represents a nontrivial class in $H_2(D_h, \partial D_h)$ because $(\pi_i \circ F)_* : H_2(\Sigma, \partial\Sigma) \to H_2(D_h, \partial D_h)$ is a nonzero map. By the Almgren isomorphism theorem, $K$ represents a nonzero element of $\pi_1(\Cyc_1(D_h, \partial D_h; \Z))$, and thus $K_* : H_1(S^1) \to H_1(\Cyc_1(D_h, \partial D_h; \Z))$ is also a nonzero map. By definition, $K$ gives a sweepout of $D_h$ by relative 1-cycles.

\begin{figure}[h]
    \centering
    \includegraphics[width=0.25\linewidth]{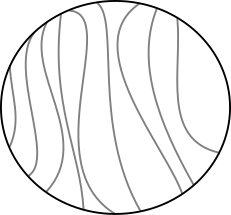}
    \qquad
    \includegraphics[width=0.25\linewidth]{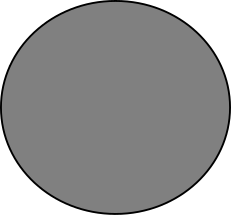}
    \caption{A 1-parameter family of 1-cycles in a Riemannian 2-disk (left) that ``glues'' into a relative 2-cycle (right), 
    which in this case represents the relative fundamental class of the disk.}
    \label{fig:GluingHomomorphism}
\end{figure}

By \cref{lem:DiskLargeWidth}, some $K(s)$ must be longer than $L+1$. As $K(s)$ is an orthogonal projection of part of $F(S^1_{st})$ onto $D_h$, we obtain that $F(S^1_{st})$ must be longer than $L + 1$ as well. This gives a contradiction.
\end{proof}

\section{Orthogonal Geodesic Chords}

\begin{figure}[h]
    \centering

    \includegraphics[width=0.7\linewidth]{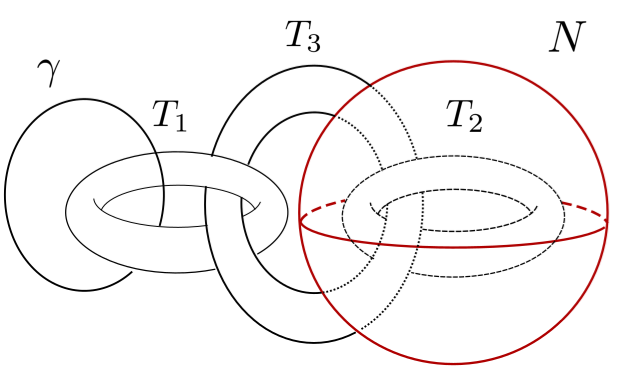}
    \caption{$N$ separates $S^3$ into two connected components, each containing one solid torus. The remaining solid torus intersects $N$ in two disks.}
    \label{fig:GeodesicsNormalToSurface}
\end{figure}

To prove \cref{thm:OrthoGeodesicChords}, we will construct a sequence of Riemannian 3-spheres $\bar{M}_h = (S^3, \bar{g}_h)$ with small diameter and volume, so that $\lambda_\rel(\bar{M}_h,N)$ and $\lambda_\rel(\bar{M}_h,\gamma)$ are large for some embedded 2-sphere $N$ and embedded circle $\gamma$. Consider three solid tori $T_1$, $T_2$, and $T_3$ embedded in $S^3$ and linked as shown in \cref{fig:GeodesicsNormalToSurface}. The embedded 2-sphere $N$ is chosen so that it intersects $T_2$ in two 2-disks. $N$ also separates $S^3$ into two closed 2-balls, $B_1$ and $B_2$, so that $B_i$ contains $T_i$ in its interior. $\gamma$ is chosen so that it links with $T_1$ as shown in the figure.
 
Similarly to our previous construction of the metric $g_h$ from \cref{sec:Construction}, give each solid torus the product metric on $D_h \times S^1$, where the length of the $S^1$ is sufficiently small so that the volume of the product metric is at most $\frac1{10}$. Extend this metric on $T_1 \cup T_2 \cup T_3$ to the entire $S^3$ as in \cref{sec:Construction}, so that it is sufficiently small away from some open neighbourhood of $T_1 \cup T_2 \cup T_3$. This defines a metric $\bar{g}_h$ on $S^3$ with diameter and volume at most 1. As a Riemannian submanifold, $N$ has diameter at most 1 due to our choice of metric on $\bar{M}_h$.

Recall from the Introduction that $\Omega_N S^3$ denotes the space of piecewise smooth paths in $S^3$ whose endpoints lie on $N$. One can prove that
\begin{equation}
    \label{eq:RelativePathsHomotopyGroups}
    \pi_i(\Omega_N S^3, \Lambda^0 N) \iso \pi_{i+1}(S^3, N),
\end{equation}
by generalizing the proof that $\pi_i(\Omega S^3) \iso \pi_{i+1}(S^3)$. Thus we have $\pi_1(\Omega_N S^3, \Lambda^0 N) = 0$ but $\pi_2(\Omega_N S^3, \Lambda^0 N) \neq 0$.

Consider some continuous map $f : (D^2, \partial D^2) \to (\Omega_N S^3, \Lambda^0 N)$ that represents a nonzero class in $\pi_2(\Omega_N S^3, \Lambda^0 N)$. Similar to the discussion in the Plan of the Proof, $f$ induces a map $F : (D^2 \times I/{\sim}) \homot D^3 \to S^3$, where $(p,r) \sim (p,r')$ for all $p \in \partial D^2$. The quotient map $q : D^2 \times I \to D^3$ can be chosen so that it sends $p \times I$ to $D^3 \cap (\R \times \{p\})$, where $D^3$ is regarded as the unit 3-disk in $\R^4$. Define $D^1_{st} = D^3 \cap (\R \times \{(s,t)\})$ and $D^2_t = D^3 \cap (\R \times \{t\})$.

Suppose that every $F(D^1_{st})$ is shorter than $L$ for some $L > 0$. By taking the double of $F$ to get a map $S^3 \to S^3$, we may apply \cref{lem:morse} to conclude that for any $\delta > 0$ we can approximate $F$ by a smooth map $\hat{F} : D^3 \to S^3$ that is transverse to $Y = \partial T_1 \cup \partial T_2$ so that:
\begin{enumerate}
    \item[I.] For each $s$ and $t$, the curves $F(D^1_{st})$ and $\hat{F}(D^1_{st})$ differ in length by less than $\delta$.

    \item[II.] Each set $\hat{F}^{-1}(Y) \cap D^2_t$ is a level set of a Morse function on $\hat{F}^{-1}(Y)$.
\end{enumerate}

In particular, each $X_i = \hat{F}^{-1}(\partial T_i)$ is a closed submanifold of $S^3$, which we endow with the preimage orientation. Henceforth we will write $F$ to mean $\hat{F}$.

\begin{lemma}
    \label{lem:OGC_DegreeOverTorus}
    For some $i = 1,2$, the map $F|_{X_i} : X_i \to \partial T_i$ has nonzero degree.
\end{lemma}
\begin{proof}
    The isomorphism $\pi_2(\Omega_N S^3, \Lambda^0 N) \iso \pi_3(S^3,N)$ from \cref{eq:RelativePathsHomotopyGroups} sends $[f] \neq 0 $ to $[F] \neq 0$. The Hurewicz homomorphism $\varphi : \pi_3(S^3, N) \to H_3(S^3, N)$ sends $[F]$ to a class $\alpha \in H_3(S^3, N) \iso \tilde{H}_3(S^3/N) \iso \Z \oplus \Z$. Statement (II) above implies that $\alpha = (k_1, k_2) \in \Z \oplus \Z$, where $k_i = \deg F|_{X_i}$. (To see this, note that the preimages of a regular value of $F$ have neighbourhoods that are foliated by the $D^2_t$'s.) Therefore it suffices to prove that $\varphi$ is an isomorphism. Consider the following commutative diagram whose rows are long exact sequences and whose columns are Hurewicz maps:
    \begin{equation}
    \begin{tikzcd}
        \pi_3(N) \dar \rar & \pi_3(S^3) \dar{\varphi_3} \rar & \pi_3(S^3, N)  \dar{\varphi} \rar & \pi_2(N)  \dar{\varphi_2} \rar & \overbrace{\pi_2(S^3)}^0  \dar
        \\
        \underbrace{H_  3(N)}_0 \rar & H_3(S^3) \rar & H_3(S^3, N) \rar & H_2(N) \rar & \underbrace{H_2(S^3)}_0
    \end{tikzcd}
    \end{equation}
    $\varphi_2$ and $\varphi_3$ are isomorphisms by the Hurewicz theorem. The Five Lemma implies that $\varphi$ is also an isomorphism. (We use a stronger version of the Five Lemma \cite[p.~129]{Hatcher_AlgTop}, in which the leftmost column is only required to be surjective and the rightmost column is only required to be injective.)
\end{proof}

\begin{proof}[Proof of \cref{thm:OrthoGeodesicChords}]
    Consider any $E > 0$. The H\"older inequality implies that $\length(\gamma) \leq \sqrt{E(\gamma)}$ for piecewise smooth curves $\gamma$ parametrized over $I$. \Cref{lem:DiskLargeWidth} allows us to choose some $h$ such that $\width_1^1(D_h) > \sqrt{E} + 1$. Consider the Riemannian 3-sphere $M = (S^3, \bar{g}_h)$ defined at the beginning of this section.
    
    Let us first prove that $\lambda_\rel(M, N) > E$. Suppose for the sake of contradiction that $\lambda_\rel(M, N) \leq E$. Then some nonzero class in $\pi_2(\Omega_N M, N)$ is represented by a map $f : (D^2, \partial D^2) \to (\Omega_N M, N)$ such that each $f(p)$ is shorter than $\sqrt{E} + \frac12$. As explained previously, \cref{lem:morse} implies that $f$ corresponds to a map $F$ satisfying statements (I) and (II), and so that every curve $F(D^1_{st})$ is shorter than $\sqrt{E} + \frac12$.

    Without loss of generality, \cref{lem:OGC_DegreeOverTorus} implies that $F|_{X_1} : X_1 \to \partial T_1$ has nonzero degree. An argument adapted from the proof of \cref{lem:winding} implies that for some $t$, $X_1 \cap D^2_t$ contains some embedded circle $C$ such that $F|_C : C \to \partial T_1$ is not nullhomotopic. Thus $F(C)$ winds around the core curve of $T_j$ a nonzero number of times, where $j = 1$ or 3. (Here we used the fact that $T_1$ and $T_3$ are linked.) Let $\pi_j : T_j \to D_h$ denote the canonical projection. The proof of \cref{lem:SurfaceExistence} works nearly verbatim to prove that $F^{-1}(T_j) \cap D^2_t$ contains a surface $\Sigma$ such that $(\pi_j \circ F)_* : H_2(\Sigma, \partial\Sigma) \to H_2(D_h, \partial D_h)$ is a nonzero map. (The crucial fact is that $C$ bounds a disk in $D^2_t$.)
    
    As in the proof of \cref{thm:theorem}, $D_h$ is swept out by a family of relative 1-cycles $\pi_j(F(D^1_{st} \cap \Sigma))$ so as a consequence of \cref{lem:DiskLargeWidth}, one of the curves $F(D^1_{st})$ must be longer than $\sqrt{E} + 1$, giving a contradiction.

    Next we prove that $\lambda_\rel(M, \gamma) > E$. Note that $\pi_1(\Omega_\gamma M, \Lambda^0 \gamma) \iso \pi_2(M, \gamma) \neq 0$, so we consider a map of pairs $h : (I, \partial I) \to (\Omega_\gamma M, \Lambda^0 \gamma)$ that represents a nonzero class in $\pi_1(\Omega_\gamma M, \Lambda^0 \gamma)$. Suppose for the sake of contradiction that every curve $h(p)$ has length at most $\sqrt{E} + \frac12$. Similar to previous arguments, $h$ induces a map $H : D^2 \homot (I \times I/{\sim}) \to M$, where $(p,r) \sim (p,r')$ for all $p \in \partial I$. The quotient map $I \times I \to D^2$ can be chosen to send $\{t\} \times I$ to $I_t = D^2 \cap (\R \times \{t\})$. $H$ represents a nonzero class in $\pi_2(M, \gamma)$, and each $H(I_t)$ has length at most $\sqrt{E} + \frac12$.
    
    The long exact sequence of homotopy groups of the pair $(M, \gamma)$ reveals that the boundary map $\pi_2(M, \gamma) \to \pi_1(\gamma)$ is an isomorphism, so $H|_{\partial D^2}$ winds around $\gamma$ a nonzero number of times. Since $\gamma$ is linked with the core curve of $T_1$, the surface $H(D^2)$ has a nonzero intersection number with that core curve. Similarly to the previous arguments, $H$ can be perturbed to a homotopic map that is transverse to $\partial T_1$ while changing the lengths of the curves $H(I_t)$ only slightly. The argument used to prove \cref{lem:SurfaceExistence} implies the existence of some surface $\Sigma \subset H^{-1}(T_1)$ such that $(\pi_1 \circ H)_* : H_2(\Sigma, \partial\Sigma) \to H_2(D_h, \partial D_h)$ is a nonzero map. As before, $D_h$ is now swept out by relative 1-cycles $\pi_j(H(I_t \cap \Sigma))$, and similar arguments as before show that one of the curves $H(I_t)$ must be longer than $\sqrt{E} + 1$, giving a contradiction.
\end{proof}

\section*{Acknowledgements}

The authors would like to thank Regina Rotman for suggesting the problem to us, and for helpful conversations. The first author was supported by the University of Toronto Excellence Award. The second author was supported by the Vanier Canada Graduate Scholarship.

\bibliography{References}

\end{document}